\documentclass[12pt]{article}

\usepackage{amssymb,amsmath} 

\DeclareMathAlphabet{\eufrak}{U}{}{}{} 
\SetMathAlphabet\eufrak{normal}{U}{euf}{m}{n}
\SetMathAlphabet\eufrak{bold}{U}{euf}{b}{n}

\numberwithin{equation}{section}

\newenvironment{Proof}{\removelastskip\par\medskip
\noindent{\em Proof.} \rm}{\penalty-20\null\hfill$\square$\par\medbreak}

 \allowdisplaybreaks

 \def\real{{\mathord{\mathbb R}}}
 
 \def\inte{{\mathord{\mathbb N}}}
 
 \def\qu{{\mathord{\mathbb Z}}}
 
 \def\Cov{{\mathrm{{\rm Cov}}}}
 
 \def\Dom{{\mathrm{{\rm Dom}}}}

 \def\real{{\mathord{{\rm I\kern-3pt R}}}}        
 
 \def\inte{{\mathord{{\rm I\kern-3pt N}}}}
 \def\sZZ{{\rm Z\kern-.45em{}Z}}
 \def\sQQ{{\kern 0.27em \vrule height1.45ex width0.03em depth0em
           \kern-0.30em \rm Q}}
 \def\qu{{\mathchoice
         {\sQQ}
         {\sQQ}
   {\kern 0.225em \vrule height1.05ex width0.025em depth0em \kern-0.25em \rm Q}
   {\kern 0.180em \vrule height0.78ex width0.020em depth0em \kern-0.20em \rm Q}
         }}
 \def\sGG{{\kern 0.27em \vrule height1.45ex width0.03em depth0em
           \kern-0.30em \rm G}}
 \def\gg{{\mathchoice
         {\sGG}
         {\sGG}
   {\kern 0.225em \vrule height1.05ex width0.025em depth0em \kern-0.25em \rm G}
   {\kern 0.180em \vrule height0.78ex width0.020em depth0em \kern-0.20em \rm G}
         }}

 \newtheorem{prop}{Proposition}[section]

 \newtheorem{corollary}[prop]{Corollary}

 \def\Dom{{\mathrm{{\rm Dom \! \ }}}}

\def\E{\mathop{\hbox{\rm I\kern-0.20em E}}\nolimits}

 \newcounter{hyp}
\title{\huge Moments of Poisson stochastic integrals with random integrands} 
 
\author{
\large Nicolas Privault
\\ 
\normalsize 
Division of Mathematical Sciences 
\\ 
\normalsize 
School of Physical and Mathematical Sciences 
\\ 
\normalsize 
Nanyang Technological University 
\\ 
\normalsize 
21 Nanyang Link 
\\ 
\normalsize 
Singapore 637371
}

\begin{document}

\hyphenation{func-tio-nals} 
\hyphenation{Privault} 

\maketitle 

\vspace{-0.9cm}

\baselineskip0.6cm
 
\begin{abstract} 
 We show that the moment of order $n$ of the Poisson 
 stochastic integral of a random 
 process $(u_x)_{x \in X}$ over a metric space $X$ 
 is given by 
\begin{equation} 
\label{pp} 
 E \left[ 
 \left( \int_X u_x(\omega ) \omega (dx) \right)^n 
 \right] 
 = 
 \sum_{ P_1,\ldots,P_k} 
 E \left[ 
 \int_{X^k} 
 \varepsilon^+_{ \eufrak{s}_k} 
 ( 
 u^{|P_1|}_{s_1} 
 \cdots 
 u^{|P_k|}_{s_k} 
 ) 
 \sigma ( d s_1 ) \cdots \sigma ( d s_k ) 
 \right] 
, 
\end{equation} 
 where the sum runs over all partitions  
 $P_1\cup \cdots \cup P_k$ of $\{ 1, \ldots , n \}$, 
 $|P_i|$ denotes the cardinality of $P_i$, 
 and $\varepsilon^+_{ \eufrak{s}_k}$ is the operator that acts 
 by addition of points at $s_1,\ldots ,s_k$ to Poisson configurations. 
 This formula 
 recovers known results in case 
 $(u(x) )_{x\in X}$ is a deterministic function on $X$. 
\end{abstract} 
 
\noindent {\bf Key words:} Poisson stochastic integrals, moment identities, Skorohod integral. 
\\ 
{\em Mathematics Subject Classification (2010):} 60G57, 60G55, 60H07. 
 
\baselineskip0.7cm

\section{Introduction} 
 Let $\Omega^X$ denote the configuration space on 
 a $\sigma$-compact metric space $X$ with 
 Borel $\sigma$-algebra ${\cal B}(X)$ and a 
 $\sigma$-finite diffuse measure $\sigma$, 
 i.e. $\Omega^X$ is 
 the space of at most countable and locally finite 
 subsets 
 of $X$, 
 defined as 
$$\Omega^X = \left\{ 
 \omega = 
 ( x_i )_{i=1,\ldots,N} \subset X, \ x_i\not= x_j \ 
 \forall i\not= j, \ N \in \inte \cup \{ \infty \} \right\} 
. 
$$ 
 Each element $\omega$ of $\Omega^X$ is identified with the 
 Radon point measure 
$$ 
 \omega = \sum_{i=1}^{\omega (X)} \epsilon_{x_i}, 
$$ 
 where $\epsilon_x$ denotes the Dirac measure at $x\in X$ 
 and $\omega (X) \in \inte \cup \{ \infty \}$ 
 represents the cardinality of $\omega (X)$. 
 The space $\Omega^X$ is endowed with the 
 Poisson probability measure $\pi_\sigma$ on $X$ 
 such that for all compact disjoint subsets $A_1,\ldots ,A_n$ of $X$, 
 $n\geq 1$, the mapping 
$$ 
 \omega \mapsto ( \omega (A_1 ) , \ldots , \omega (A_n)) 
$$ 
 is a vector of independent Poisson distributed random variables 
 on $\inte$ with respective intensities $\sigma (A_1), \ldots , \sigma (A_n)$. 
\\ 
 
 In \cite{bassan} the moment formula 
\begin{equation} 
\label{**} 
 E \left[ 
 \left( \int_X f (x) \omega (dx) \right)^n 
 \right] 
 = 
 n! 
 \sum_{ 
 r_1+2r_2+\cdots + n r_n =n 
 \atop 
 r_1,\ldots , r_n \geq 0
} 
 \prod_{k=1}^n \frac{1}{(k!)^{r_k}r_k!} 
 \prod_{k=1}^n 
 \left( \int_X f^k(x) \sigma (dx) \right)^{r_k} 
\end{equation} 
 has been proved for $f:X\to \real$ a deterministic 
 sufficiently integrable function. 
 The proof of \cite{bassan} relies on the L\'evy-Khintchine 
 representation of the Laplace transform of 
 $\int_X u(x) \omega (dx)$, and this result can also be recovered 
 under a different combinatorial interpretation 
 by the Fa\`a di Bruno formula, cf. e.g. \S~2.4 of \cite{lukacs}, 
 from the relation 
\begin{equation} 
\label{ftr} 
 E \left[ 
 \left( \int_X u(x) \omega (dx) \right)^n 
 \right] 
 = 
 \sum_{ P_1,\ldots,P_k} 
 \int_X 
 u^{|P_1|} (x) 
 \sigma ( dx ) 
 \cdots 
 \int_X 
 u^{|P_k|} (x) 
 \sigma ( dx ) 
, 
\end{equation} 
 between the moments and the cumulants $\kappa_n = \int_X u^n (x) \sigma (dx)$, 
 $n\geq 1$, of $\displaystyle \int_X u(x) \omega (dx)$, 
 where the sum runs over all partitions  
 $P_1\cup \cdots \cup P_k$ of $\{ 1, \ldots , n \}$. 
 \\ 
 
 Recently, \eqref{**} has been applied to control 
 the $p$-variation and the number of crossings 
 of fractional Poisson and shot noise processes with deterministic kernels, 
 cf. \cite{bierme2}. 
\\ 
 
 In this paper we extend the above formula \eqref{**} 
 to random integrands. 
 Namely, we state that given $u:\Omega^X \times X \to \real$ 
 a sufficiently integrable random process we have 
\begin{equation} 
\label{asfg} 
 E \left[ 
 \left( \int_X u_x (\omega ) \omega (dx) \right)^n 
 \right] 
 = 
 \sum_{ P_1,\ldots,P_k} 
 E \left[ 
 \int_{X^k} 
 \varepsilon^+_{ \eufrak{s}_k} 
 \left( 
 u^{|P_1|}_{s_1} 
 \cdots 
 u^{|P_k|}_{s_k} 
 \right) 
 \sigma ( ds_1 ) \cdots \sigma ( ds_k ) 
 \right] 
, 
\end{equation} 
 cf. Proposition~\ref{p111} below, 
 where the sum runs over all (disjoint) partitions  
 $P_1\cup \cdots \cup P_k$ of $\{ 1, \ldots , n \}$, 
 $k=1,\ldots ,n$, 
 and $|P|$ denotes the cardinal of $P \subset \{1,\ldots , n \}$. 
 In \eqref{asfg}, 
 $\varepsilon^+_{s_k}$ is the addition operator 
 defined 
 on any random variable $F:\Omega^X \to \real$ 
 by 
$$ 
 \varepsilon^+_{\eufrak{s}_k} 
 F(\omega ) = F(\omega \cup \{ s_1,\ldots ,s_k \} ), 
 \qquad 
 \omega \in \Omega^X, 
 \quad 
 s_1,\ldots , s_k \in X, 
$$ 
 where 
$$ 
 \eufrak{s}_k = (s_1,\ldots ,s_k ) \in X^k, 
 \quad 
 k \geq 1. 
$$ 
 As expected, $(u(x) )_{x \in X}$ is a deterministic function 
 we have 
$$ 
 \varepsilon^+_{\eufrak{s}_k } 
 u(s_i) 
 = 
 u(s_i) 
, 
\qquad 
 1 \leq i \leq k, 
$$ 
 in which case \eqref{asfg} recovers \eqref{ftr}. 
\subsubsection*{Examples} 
 In the case of second order moments, 
 \eqref{asfg} yields 
$$ 
 E \left[ 
 \left( \int_X u_x (\omega ) \omega (dx) \right)^2 
 \right] 
 = 
 E \left[ 
 \int_X 
 \varepsilon^+_{ \eufrak{s}} 
 | u |^2_s 
 \sigma ( ds ) 
 \right] 
 + 
 E \left[ 
 \int_{X^2} 
 \varepsilon^+_{ \eufrak{s}_1} \varepsilon^+_{ \eufrak{s}_2} 
 ( u_{s_1} u_{s_2} ) 
 \sigma ( ds_1 ) \sigma ( ds_2 ) 
 \right] 
, 
$$ 
 which, if $\lambda : = \sigma (X) < \infty$, recovers 
\begin{eqnarray*} 
 E \left[ 
 ( \omega (X) )^6 
 \right] 
 & = & 
 E \left[ 
 \left( \int_X u_x (\omega ) \omega (dx) \right)^3 
 \right] 
\\ 
 & = & 
 E \left[ 
 \int_X 
 \varepsilon^+_{ \eufrak{s}} 
 ( \omega (X) )^2 
 \sigma ( ds ) 
 \right] 
 + 
 E \left[ 
 \int_{X^2} 
 \varepsilon^+_{ \eufrak{s}_1} \varepsilon^+_{ \eufrak{s}_2} 
 ( \omega (X) )^2 
 \sigma ( ds_1 ) \sigma ( ds_2 ) 
 \right] 
\\ 
 & = & 
 \lambda 
 E \left[ 
 ( \omega (X) + 1 )^2 
 \right] 
 + 
 \lambda^2 
 E \left[ 
 ( \omega (X) + 2 )^2 
 \right] 
\\ 
 & = & 
 \lambda 
 + 7 \lambda^2 
 + 6 \lambda^3 
 + \lambda^4 
\\ 
 & = & B_4 ( \lambda ), 
\end{eqnarray*} 
 by taking $u_x(\omega ) = \omega (X)$, 
 $x\in X$, where $B_4$ is the Bell polynomial of order 
 $4$, cf. \eqref{bn} below. 
\\ 
 
 Concerning third order moments, \eqref{asfg} shows that 
\begin{eqnarray*} 
 E \left[ 
 \left( \int_X u_x (\omega ) \omega (dx) \right)^3 
 \right] 
 & = & 
 E \left[ 
 \int_X 
 \varepsilon^+_{ \eufrak{s}} 
 | u |^3_s 
 \sigma ( ds ) 
 \right] 
 + 
 3 E \left[ 
 \int_{X^2} 
 \varepsilon^+_{ \eufrak{s}_1} \varepsilon^+_{ \eufrak{s}_2} 
 ( | u_{s_1} |^2 u_{s_2} ) 
 \sigma ( ds_1 ) \sigma ( ds_2 ) 
 \right] 
\\ 
 & & 
 + 
 E \left[ 
 \int_{X^3} 
 \varepsilon^+_{ \eufrak{s}_1} \varepsilon^+_{ \eufrak{s}_2} 
 \varepsilon^+_{ \eufrak{s}_3} 
 ( u_{s_1} u_{s_2} u_{s_3} ) 
 \sigma ( ds_1 ) \sigma ( ds_2 ) \sigma ( ds_3 ) 
 \right] 
, 
\end{eqnarray*} 
 which, if $\lambda : = \sigma (X) < \infty$, 
 and taking $u_x(\omega ) = \omega (X)$, 
 $x\in X$, yields 
\begin{eqnarray*} 
\lefteqn{ 
 E \left[ 
 ( \omega (X) )^6 
 \right] 
 = 
 E \left[ 
 \left( \int_X u_x (\omega ) \omega (dx) \right)^3 
 \right] 
} 
\\ 
 & = & 
 E \left[ 
 \int_X 
 \varepsilon^+_{ \eufrak{s}} 
 ( \omega (X) )^3 
 \sigma ( ds ) 
 \right] 
 + 
 3 
 E \left[ 
 \int_{X^2} 
 \varepsilon^+_{ \eufrak{s}_1} \varepsilon^+_{ \eufrak{s}_2} 
 ( ( \omega (X) )^2 \omega (X) ) 
 \sigma ( ds_1 ) \sigma ( ds_2 ) 
 \right] 
\\ 
 & & 
 + 
 E \left[ 
 \int_{X^3} 
 \varepsilon^+_{ \eufrak{s}_1} \varepsilon^+_{ \eufrak{s}_2}  
 \varepsilon^+_{ \eufrak{s}_3} 
 ( \omega (X) )^3 
 \sigma ( ds_1 ) \sigma ( ds_2 ) \sigma ( ds_3 ) 
 \right] 
\\ 
 & = & 
 \lambda 
 E \left[ 
 ( \omega (X) + 1 )^3 
 \right] 
 + 
 3 
 \lambda^2 
 E \left[ 
 ( \omega (X) + 2 )^2 ( \omega (X) + 2 ) 
 \right] 
\\ 
 & & 
 + 
 \lambda^3 
 E \left[ 
 ( \omega (X) + 3 )^3 
 \right] 
\\ 
 & = & 
 \lambda 
 + 31 \lambda^2 
 + 90 \lambda^3 
 + 65 \lambda^4 
 + 15 \lambda^5 
 + \lambda^6 
\\ 
 & = & B_6 ( \lambda ), 
\end{eqnarray*} 
 where $B_6$ is the Bell polynomial of order $6$. 
\\ 
 
 We proceed as follows. 
 In Section~\ref{ci} we rewrite a result of \cite{prinv} 
 on the moments of compensated Poisson-Skorohod integrals 
 in the language of set partitions. 
 In Section~\ref{nci} we deduce formulas for non-compensated 
 integrals of random integrands by binomial inversion. 
 In the case of deterministic integrands 
 and indicator functions, in Section~\ref{ifp} we recover 
 and extend known relations between the moments of the Poisson distribution 
 and Stirling numbers. 
\\ 
 
 The moment identities in this paper are stated for bounded 
 random variables $F$ and processes $u$ with compact support, 
 however 
 they can be extended by assuming suitable conditions ensuring 
 that the right hand side of the formula is finite. 
\section{Poisson-Skorohod integrals} 
\label{ci} 
 Our proof of moment identities relies 
 on the Skorohod integral operator $\delta$ 
 which is defined on any measurable 
 process $u : \Omega^X \times X \to \real$ 
 by the expression 
\begin{equation} 
\label{ui} 
 \delta ( u ) 
 = 
 \int_X u_x ( \omega \setminus x ) 
 \omega ( dx ) 
 - 
 \int_X u_x ( \omega ) \sigma (dx) 
, 
\end{equation} 
 provided 
$ 
 \displaystyle 
 E \left[ \int_X | u_x ( \omega ) | \sigma (dx ) \right] < \infty 
$, 
 cf. Corollary~1 of \cite{picard}. 
 In \eqref{ui}, 
 $\omega \setminus x$ denotes the configuration $\omega \in \Omega^X$ 
 after removal of the point $x$ in case $x\in \omega$. 
\\ 
 
 We start with a moment identity for compensated 
 Poisson-Skorohod integrals, obtained by rewriting 
 Theorem~5.1 of \cite{prinv} using set partitions. 
 By saying that 
 $u:\Omega^X \times X \to \real$ has a compact support in $X$ we 
 mean that there exists a compact subset $K$ of $X$ such that 
 $u_x (\omega ) = 0$ for all $\omega \in \Omega^X$ and $x \in X \setminus K$. 
\begin{prop} 
\label{cc} 
 Let $F:\Omega^X \to \real$ be a bounded random variable and 
 let $u:\Omega^X \times X \to \real$ be a bounded process 
 with compact support in $X$. 
 For all $n\geq 0$ we have 
$$ 
 E \left[ 
 \delta ( u )^n 
 F 
 \right] 
 = 
 \sum_{c=0}^n 
 (-1)^c 
 {n \choose c} 
 \sum_{k =0}^{n-c} 
 \sum_{ 
 \substack{ 
 l_1 + \cdots + l_k = n - c 
 \\ 
 l_1,\ldots ,l_k \geq 1 
 \\ 
 l_{k+1},\ldots ,l_{k+c} = 1 
 } 
 } 
 \! \! \! \! \! \! \! \! \! 
 {\cal N}_{\eufrak{L}_k} 
 E \left[  \int_{X^{k+c}} 
 \! \! \! 
 \varepsilon^+_{\eufrak{s}_k } 
 F 
 \prod_{p=1}^{k+c}
 \varepsilon^+_{ \eufrak{s}_k \setminus s_p} 
 u^{l_p}_{s_p} 
 \hskip0.05cm 
 d\sigma^{k+c} ( \eufrak{s}_{k+c} ) 
 \right] 
, 
$$ 
 where $d\sigma^b ( \eufrak{s}_b ) = \sigma ( ds_1 ) \cdots \sigma ( ds_b )$, 
 $\eufrak{L}_k = (l_1,\ldots ,l_k )$, and ${\cal N}_{\eufrak{L}_k}$ is 
 the number of partitions of a set of $l_1+\cdots + l_k$ elements 
 into $k$ subsets of lengths $l_1,\ldots ,l_k \geq 1$. 
\end{prop} 
\begin{Proof} 
 The proof of this formula relies on the identity 
\begin{equation} 
\label{plkmn} 
 E \left[ 
 \delta ( u )^n 
 F 
 \right] 
 = 
 \sum_{k=0}^n 
 \sum_{b = k}^n 
 (-1)^{b-k} 
 \! \! \! \! \! \! \! \! \! 
 \sum_{ 
 \substack{ 
 l_1 + \cdots + l_k = n - (b-k) 
 \\ 
 l_1,\ldots ,l_k \geq 1 
 \\ 
 l_{k+1},\ldots ,l_b = 1 
 } 
 } 
 \! \! \! \! \! \! \! \! \! 
 C_{\eufrak{L}_k,b} 
 E \left[  \int_{X^b} 
 \! \! \! 
 \varepsilon^+_{\eufrak{s}_k } 
 F 
 \prod_{p=1}^b 
 \varepsilon^+_{ \eufrak{s}_k \setminus s_p} 
 u^{l_p}_{s_p} 
 \hskip0.05cm 
 d\sigma^b ( \eufrak{s}_b ) 
 \right] 
\end{equation} 
 for the moments of the 
 compensated Poisson-Skorohod integral $\delta (u)$, cf. 
 Theorem~5.1 of \cite{prinv} and Theorem~1 of \cite{priinvcr}, 
 where 
\begin{equation} 
\label{defin} 
 C_{\eufrak{L}_k,k+c} 
 = 
 \! \! \! \! \! \! \! 
 \sum_{0 = r_{c+1} < \cdots < r_0 = k+c+1} 
 \ 
 \prod_{q=0}^c 
 \ 
 \prod_{p=r_{q+1} + q - c + 1 }^{r_q+q-c-1} 
 {l_1 + \cdots + l_p + q - 1  
 \choose l_1 + \cdots + l_{p-1} + q } 
. 
\end{equation} 
 Next we note that 
 $C_{\eufrak{L}_k,k+c}$ defined in \eqref{defin} represents 
 the number of partitions of a set of $n=l_1+\cdots + l_k+c$ elements 
 into $k$ subsets of lengths $l_1,\ldots ,l_k$ and $c$ singletons, 
 hence when $l_1 + \cdots + l_k = n - c$ we have 
$$ 
 C_{\eufrak{L}_k,k+c} 
 = 
 {n \choose c} 
 {\cal N}_{\eufrak{L}_k} 
, 
$$ 
 since $ {\cal N}_{\eufrak{L}_k}$ is 
 the number of partitions of a set of $l_1+\cdots + l_k=n-c$ elements 
 into $k$ subsets of lengths $l_1,\ldots ,l_k$. 
 Hence we have 
\begin{eqnarray*} 
\lefteqn{ 
 E \left[ 
 \delta ( u )^n 
 F 
 \right] 
 = 
 \sum_{k=0}^n 
 \sum_{b = k}^n 
 (-1)^{b-k} 
 \! \! \! \! \! \! \! \! \! 
 \sum_{ 
 \substack{ 
 l_1 + \cdots + l_k = n - ( b - k ) 
 \\ 
 l_1,\ldots ,l_k \geq 1 
 \\ 
 l_{k+1},\ldots ,l_b = 1 
 } 
 } 
 \! \! \! \! \! \! \! \! \! 
 C_{\eufrak{L}_k,b} 
 E \left[  \int_{X^b} 
 \! \! \! 
 \varepsilon^+_{\eufrak{s}_k } 
 F 
 \prod_{p=1}^b 
 \varepsilon^+_{ \eufrak{s}_k \setminus s_p} 
 u^{l_p}_{s_p} 
 \hskip0.05cm 
 d\sigma^b ( \eufrak{s}_b ) 
 \right] 
} 
\\ 
 & = & 
 \sum_{c=0}^n 
 \sum_{k =0}^{n-c} 
 (-1)^c 
 \! \! \! \! \! \! \! \! \! 
 \sum_{ 
 \substack{ 
 l_1 + \cdots + l_k = n - c 
 \\ 
 l_1,\ldots ,l_k \geq 1 
 \\ 
 l_{k+1},\ldots ,l_{k+c} = 1 
 } 
 } 
 \! \! \! \! \! \! \! \! \! 
 C_{\eufrak{L}_k,k+c} 
 E \left[  \int_{X^{k+c}} 
 \! \! \! 
 \varepsilon^+_{\eufrak{s}_k } 
 F 
 \prod_{p=1}^{k+c}
 \varepsilon^+_{ \eufrak{s}_k \setminus s_p} 
 u^{l_p}_{s_p} 
 \hskip0.05cm 
 d\sigma^{k+c} ( \eufrak{s}_{k+c} ) 
 \right] 
\\ 
 & = & 
 \sum_{c=0}^n 
 (-1)^c 
 {n \choose c} 
 \sum_{k =0}^{n-c} 
 \sum_{ 
 \substack{ 
 l_1 + \cdots + l_k = n - c 
 \\ 
 l_1,\ldots ,l_k \geq 1 
 \\ 
 l_{k+1},\ldots ,l_{k+c} = 1 
 } 
 } 
 \! \! \! \! \! \! \! \! \! 
 {\cal N}_{\eufrak{L}_k} 
 E \left[  \int_{X^{k+c}} 
 \! \! \! 
 \varepsilon^+_{\eufrak{s}_k } 
 F 
 \prod_{p=1}^{k+c}
 \varepsilon^+_{ \eufrak{s}_k \setminus s_p} 
 u^{l_p}_{s_p} 
 \hskip0.05cm 
 d\sigma^{k+c} ( \eufrak{s}_{k+c} ) 
 \right] 
. 
\end{eqnarray*} 
\end{Proof} 
 The proof of \eqref{plkmn} relies on the duality relation 
\begin{equation} 
\label{dr1} 
 E [ \langle D F , u \rangle_{L^2 ( X )}] 
 = 
 E [ F \delta ( u ) ], 
\end{equation} 
 between $\delta$ and the 
 finite difference gradient 
\begin{equation} 
\label{prr} 
 D_x F(\omega ) = \varepsilon^+_x F (\omega ) - F( \omega ), 
 \qquad 
 \omega \in \Omega^X, 
 \quad 
 x\in X, 
\end{equation} 
 for all $F$ and $u$ in the respective 
 closed $L^2$ domains 
 $\Dom ( \delta ) 
 \subset L^2 ( \Omega^X \times X , \pi_\sigma \otimes \sigma)$ 
 and $\Dom ( D ) \subset L^2 ( \Omega^X , \pi_\sigma )$ 
 of $D$ and $\delta$, cf. \cite{prinv} and references 
 therein. 
\section{Pathwise integrals} 
\label{nci} 
 The next Proposition~\ref{p111} is the main result of this paper, 
 it yields \eqref{asfg} and 
 follows directly by binomial inversion of 
 Proposition~\ref{cc}. We state its proof 
 due to the additional presence of expectations. 
\begin{prop} 
\label{p111} 
 Let $F:\Omega^X \to \real$ be a bounded random variable, and let 
 $u:\Omega^X \times X \to \real$ be a bounded random process 
 with compact support in $X$. For all $n\geq 0$ we have 
$$ 
 E \left[ 
 F 
 \left( \int_X u_x (\omega ) \omega (dx) \right)^n 
 \right] 
 = 
 \sum_{ P_1,\ldots,P_k} 
 E\left[ 
 \int_{X^k} 
 \varepsilon^+_{ \eufrak{s}_k } 
 ( 
 F 
 u^{|P_1|}_{s_1} 
 \cdots 
 u^{|P_k|}_{s_k} 
 ) 
 \hskip0.05cm 
 \sigma ( ds_1 ) 
 \cdots 
 \sigma ( ds_k ) 
 \right] 
. 
$$ 
\end{prop} 
\begin{Proof} 
 We have, applying Proposition~\ref{cc} at the rank $n-i$ 
 to the process 
 $\varepsilon^+ u = (\varepsilon^+_x u_x )_{x\in X}$, 
\begin{eqnarray} 
\nonumber 
\lefteqn{ 
 E \left[ 
 F 
 \left( \int_X u_x ( \omega ) \omega (dx) \right)^n 
 \right] 
 = 
 E \left[ 
 F 
 \left( \delta ( \varepsilon^+ u ) + \int_X \varepsilon^+_x u_x 
 \sigma (dx) \right)^n 
 \right] 
} 
\\ 
\nonumber 
 & = & 
 \sum_{i=0}^n 
 {n \choose i} 
 E \left[ 
 F 
 ( \delta ( \varepsilon^+ u ) )^{n-i} 
 \left( \int_X \varepsilon^+_x u_x \sigma (dx) \right)^i 
 \right] 
\\ 
\nonumber 
 & = & 
 \sum_{i=0}^n 
 {n \choose i} 
 \sum_{c=0}^{n-i} 
 (-1)^c 
 {{n-i} \choose c} 
 \sum_{k =0}^{{n-i}-c} 
\\ 
\nonumber 
 & & 
 \! \! \! \! \! \! \! \! \! \! \! 
 \sum_{ 
 \substack{ 
 l_1 + \cdots + l_k = {n-i} - c 
 \\ 
 l_1,\ldots ,l_k \geq 1 
 } 
 } 
 \! \! \! \! \! \! \! \! \! 
 {\cal N}_{\eufrak{L}_k} 
 E \left[ 
 \int_{X^k} 
 \varepsilon^+_{\eufrak{s}_k } 
 \left( 
 F \left( \int_X \varepsilon^+_x u_x \sigma (dx) \right)^i 
 \left( 
 \int_X \varepsilon^+_x u_x \sigma (dx) 
 \right)^c 
 \right) 
 \prod_{p=1}^k 
 \varepsilon^+_{ \eufrak{s}_k \setminus s_p} 
 \varepsilon^+_{s_p} u_{s_p}^{l_p} 
 \hskip0.05cm 
 d\sigma^k ( \eufrak{s}_k ) 
 \right] 
\\ 
\nonumber 
 & = & 
 \sum_{a=0}^n 
 \sum_{i=0}^a 
 {n \choose i} 
 (-1)^{a-i} 
 {{n-i} \choose a-i} 
 \sum_{k =0}^{{n}-a} 
\\ 
\nonumber 
 & & 
 \! \! \! \! \! \! \! \! \! \! \! 
 \sum_{ 
 \substack{ 
 l_1 + \cdots + l_k = {n} - a 
 \\ 
 l_1,\ldots ,l_k \geq 1 
 } 
 } 
 \! \! \! \! \! \! \! \! \! 
 {\cal N}_{\eufrak{L}_k} 
 E \left[ 
 \int_{X^k} 
 \varepsilon^+_{\eufrak{s}_k } 
 \left( 
 F \left( \int_X \varepsilon^+_x u_x \sigma (dx) \right)^i 
 \left( 
 \int_X \varepsilon^+_x u_x \sigma (dx) 
 \right)^{a-i} 
 \right) 
 \prod_{p=1}^k 
 \varepsilon^+_{ \eufrak{s}_k } 
 u_{s_p}^{l_p} 
 \hskip0.05cm 
 d\sigma^k ( \eufrak{s}_k ) 
 \right] 
\\ 
\nonumber 
 & = & 
 \sum_{a=0}^n 
 {n \choose a} 
 \sum_{k =0}^{{n}-a} 
 \sum_{i=0}^a 
 (-1)^{a-i} 
 {a \choose i} 
\\ 
\nonumber 
 & & 
 \sum_{ 
 \substack{ 
 l_1 + \cdots + l_k = {n} - a 
 \\ 
 l_1,\ldots ,l_k \geq 1 
 } 
 } 
 \! \! \! \! \! \! \! \! \! 
 {\cal N}_{\eufrak{L}_k} 
 E \left[ 
 \int_{X^k} 
 \varepsilon^+_{\eufrak{s}_k } 
 \left( 
 F \left( \int_X \varepsilon^+_x u_x \sigma (dx) \right)^a 
 \right) 
 \varepsilon^+_{ \eufrak{s}_k } 
 \prod_{p=1}^k 
 u^{l_p}_{s_p} 
 \hskip0.05cm 
 d\sigma^k ( \eufrak{s}_k ) 
 \right] 
\\ 
\label{ascf} 
 & = & 
 \sum_{k =0}^n 
 \sum_{ 
 \substack{ 
 l_1 + \cdots + l_k = n 
 \\ 
 l_1,\ldots ,l_k \geq 1 
 } 
 } 
 {\cal N}_{\eufrak{L}_k} 
 E \left[ 
 \int_{X^k} 
 \varepsilon^+_{ \eufrak{s}_k } 
 \left( 
 F 
 \prod_{p=1}^k 
 u^{l_p}_{s_p} 
 \right) 
 \hskip0.05cm 
 d\sigma^k ( \eufrak{s}_k ) 
 \right] 
, 
\end{eqnarray}  
 since 
 $\displaystyle  \sum_{i=0}^a (-1)^{a-i} {a \choose i} = {\bf 1}_{\{ a=0 \}}$ 
 with the convention $0^0=1$. 
\end{Proof} 
 When $u : X \to \real$ is a deterministic function, 
 Proposition~\ref{cc} yields 
\begin{eqnarray*} 
 E \left[ 
 F 
 \left( \int_X u(x) \omega (dx) \right)^n 
 \right] 
 & = & 
 \sum_{ 
 P_1,\ldots,P_a} 
 \int_{X^a} 
 u^{|P_1|}_{s_1} \cdots u^{|P_a|}_{s_a}
 E\left[ 
 \varepsilon^+_{\eufrak{s}_a } 
 F 
 \right] 
 \sigma ( ds_1 ) \cdots \sigma ( ds_a ) 
, 
\end{eqnarray*} 
 which recovers \eqref{**} by taking $F=1$, 
 and 
$$ 
 \Cov \left( 
 F 
 , 
 \left( \int_X u(x) \omega (dx) \right)^m 
 \right) 
 = 
 \sum_{ 
 P_1,\ldots,P_a} 
 \int_{X^a} 
 u^{|P_1|}_{s_1} \cdots u^{|P_a|}_{s_a} 
 E\left[ 
 ( \varepsilon^+_{\eufrak{s}_a } 
 F - F ) 
 \right] 
 \sigma ( ds_1 ) \cdots \sigma ( ds_a ) 
. 
$$ 
 By \eqref{ascf}, Proposition~\ref{p111} also rewrites for 
 compensated integrals as follows. 
\begin{prop} 
 Let $F:\Omega^X \to \real$ be a bounded random variable, and let 
 $u:\Omega^X \times X \to \real$ be a bounded random process 
 with compact support in $X$. For all $n\geq 0$ we have 
\begin{eqnarray} 
\label{note} 
\lefteqn{ 
 E \left[ 
 F 
 \left( \int_X u_x ( \omega ) ( \omega (dx) - \sigma (dx ) ) \right)^n 
 \right] 
 = 
 \sum_{c=0}^n 
 (-1)^c 
 {n \choose c} 
} 
\\ 
\nonumber 
 & & 
 \! \! \! \! \! \! \! \! \! \! \! \! \! 
 \sum_{ P_1,\ldots,P_a \subset \{ 1 , \ldots , n - c \} } 
 E\left[ 
 \int_{X^a} 
 \varepsilon^+_{\eufrak{s}_a } 
 \left( 
 F 
 \left( 
 \int_X u_x ( \omega ) \sigma (dx) 
 \right)^c 
 u^{|P_1|}_{s_1} 
 \cdots 
 u^{|P_a|}_{s_a} 
 \right) 
 \hskip0.05cm 
 d\sigma^a ( \eufrak{s}_a ) 
 \right] 
. 
\end{eqnarray} 
\end{prop} 
\begin{Proof} 
 By \eqref{ascf} we have 
\begin{eqnarray*} 
\nonumber 
\lefteqn{ 
 E \left[ 
 F 
 \left( \int_X u_x ( \omega ) ( \omega (dx) - \sigma (dx ) ) \right)^n 
 \right] 
} 
\\ 
\nonumber 
 & = & 
 \sum_{c=0}^n 
 (-1)^c 
 {n \choose c} 
 E \left[ 
 F 
 \left( \int_X u_x ( \omega ) \omega (dx) ) \right)^{n-c} 
 \left( \int_X u_x ( \omega ) \sigma (dx) ) \right)^c 
 \right] 
\\ 
\nonumber 
 & = & 
 \sum_{c=0}^n 
 (-1)^c 
 {n \choose c} 
 \sum_{a =0}^{n-c} 
 \sum_{ 
 \substack{ 
 l_1 + \cdots + l_a = n - c 
 \\ 
 l_1,\ldots ,l_a \geq 1 
 \\ 
 l_{a+1},\ldots ,l_{a+c} = 1 
 } 
 } 
 \! \! \! \! \! \! \! \! \! 
 {\cal N}_{\eufrak{L}_a} 
 E \left[  \int_{X^{a+c}} 
 \! \! \! 
 \varepsilon^+_{\eufrak{s}_a } 
 \left( 
 F 
 \prod_{p=1}^{a+c}
 u^{l_p}_{s_p} 
 \right) 
 \hskip0.05cm 
 d\sigma^{a+c} ( \eufrak{s}_{a+c} ) 
 \right] 
\\ 
\nonumber 
 & = & 
 \sum_{c=0}^n 
 (-1)^c 
 {n \choose c} 
\\ 
 & & 
 \! \! \! \! \! \! \! \! \! \! \! \! \! 
 \sum_{ P_1,\ldots,P_a \subset \{ 1 , \ldots , n - c \} } 
 E\left[ 
 \int_{X^a} 
 \varepsilon^+_{\eufrak{s}_a } 
 \left( 
 F 
 \left( 
 \int_X u_x ( \omega ) \sigma (dx) 
 \right)^c 
 u^{|P_1|}_{s_1} 
 \cdots 
 u^{|P_a|}_{s_a} 
 \right) 
 \hskip0.05cm 
 d\sigma^a ( \eufrak{s}_a ) 
 \right] 
. 
\end{eqnarray*} 
\end{Proof} 
 The next proposition specializes the above result 
 to the case of deterministic integrands. 
 We note that it can also be obtained independently 
 as in \eqref{ftr} from the relation 
 between the moments and the cumulants 
 $\kappa_1=0$, $\kappa_n = \int_X u^n (x) \sigma (dx)$, 
 $n\geq 2$, 
 of $\displaystyle \int_X u(x) ( \omega (dx) - \sigma (dx))$. 
\begin{prop} 
\label{p1111} 
 Let $f : X \to \real$ be a bounded deterministic function 
 with compact support on $X$. 
 For all $n\geq 1$ we have 
$$ 
 E \left[ 
 \left( \int_X f (x) ( \omega (dx) - \sigma (dx) ) 
 \right)^n 
 \right] 
 = 
 \!  \!  \!  \!  \! 
 \sum_{ P_1,\ldots,P_k\atop 
 |P_1|\geq 2, \ldots ,|P_k|\geq 2} 
 \!  \!  \!  \!  \! 
 \int_X 
 f^{|P_1|}(x) 
 \sigma ( dx ) 
 \cdots 
 \int_X 
 f^{|P_k|}(x) 
 \sigma ( dx ) 
, 
$$ 
 where the sum runs over all partitions  
 $P_1\cup \cdots \cup P_k$ of $\{ 1, \ldots , n \}$ 
 of size at least $2$. 
\end{prop} 
\begin{Proof} 
 By Proposition~\ref{p111} and binomial inversion we have 
\begin{eqnarray*} 
\nonumber 
\lefteqn{ 
 E \left[ 
 \left( \int_X f (x) \omega (dx) - \int_X f (x) 
 \sigma (dx) \right)^n 
 \right] 
} 
\\ 
\nonumber 
 & = & 
 \sum_{c=0}^n 
 (-1)^c 
 {n \choose c} 
 \left( \int_X f (x) \sigma (dx) \right)^c 
 E \left[ 
 \left( \int_X f (x) \omega (dx) \right)^{n-c} 
 \right] 
\\ 
\nonumber 
 & = & 
 \sum_{c=0}^n 
 (-1)^c 
 {n \choose c} 
 \left( \int_X f (x) \sigma (dx) \right)^c 
 \sum_{ 
 \substack{ 
 l_1 + \cdots + l_a = n-c 
 \\ 
 l_1,\ldots ,l_a \geq 1 
 } 
 } 
 {\cal N}_{\eufrak{L}_a} 
 \int_{X^a} 
 \prod_{p=1}^a 
 f^{l_p} ({s_p}) 
 \hskip0.05cm 
 d\sigma^a ( \eufrak{s}_a ) 
\\ 
\nonumber 
 & = & 
 \sum_{c=0}^n 
 (-1)^c 
 {n \choose c} 
 \sum_{k=0}^{n-c} 
 {n-c \choose k} 
 \left( \int_X f (x) \sigma (dx) \right)^{k+c} 
 \sum_{ 
 \substack{ 
 l_1 + \cdots + l_a = n-c-k 
 \\ 
 l_1,\ldots ,l_a \geq 2 
 } 
 } 
 {\cal N}_{\eufrak{L}_a} 
 \int_{X^a} 
 \prod_{p=1}^a 
 f^{l_p} ({s_p}) 
 \hskip0.05cm 
 d\sigma^a ( \eufrak{s}_a ) 
\\ 
\nonumber 
 & = & 
 \sum_{b=0}^n 
 \sum_{c=0}^b 
 (-1)^c 
 {n \choose c} 
 {n-c \choose b-c} 
 \left( \int_X f (x) \sigma (dx) \right)^b 
 \sum_{ 
 \substack{ 
 l_1 + \cdots + l_a = n-b
 \\ 
 l_1,\ldots ,l_a \geq 2 
 } 
 } 
 {\cal N}_{\eufrak{L}_a} 
 \int_{X^a} 
 \prod_{p=1}^a 
 f^{l_p} ({s_p}) 
 \hskip0.05cm 
 d\sigma^a ( \eufrak{s}_a ) 
\\ 
\nonumber 
 & = & 
 \sum_{b=0}^n 
 {n \choose b} 
 \sum_{c=0}^b 
 (-1)^c 
 {b \choose c} 
 \left( \int_X f (x) \sigma (dx) \right)^b 
 \sum_{ 
 \substack{ 
 l_1 + \cdots + l_a = n-b
 \\ 
 l_1,\ldots ,l_a \geq 2 
 } 
 } 
 {\cal N}_{\eufrak{L}_a} 
 \int_{X^a} 
 \prod_{p=1}^a 
 f^{l_p} ({s_p}) 
 \hskip0.05cm 
 d\sigma^a ( \eufrak{s}_a ) 
\\ 
\nonumber 
 & = & 
 \sum_{b=0}^n 
 {n \choose b} 
 \sum_{c=0}^b 
 (-1)^c 
 {b \choose c} 
 \left( \int_X f (x) \sigma (dx) \right)^b 
 \sum_{ 
 \substack{ 
 l_1 + \cdots + l_a = n-b
 \\ 
 l_1,\ldots ,l_a \geq 2 
 } 
 } 
 {\cal N}_{\eufrak{L}_a} 
 \int_{X^a} 
 \prod_{p=1}^a 
 f^{l_p} ({s_p}) 
 \hskip0.05cm 
 d\sigma^a ( \eufrak{s}_a ) 
\\ 
\nonumber 
 & = & 
 \sum_{ 
 \substack{ 
 l_1 + \cdots + l_a = n
 \\ 
 l_1,\ldots ,l_a \geq 2 
 } 
 } 
 {\cal N}_{\eufrak{L}_a} 
 \int_{X^a} 
 \prod_{p=1}^a 
 f^{l_p} ({s_p}) 
 \hskip0.05cm 
 d\sigma^a ( \eufrak{s}_a ). 
\end{eqnarray*}  
\end{Proof} 
\section{Indicator functions and polynomials} 
\label{ifp} 
 When $u (x) = {\bf 1}_A (x)$ is a deterministic indicator function 
 with $A\in {\cal B}(X)$, 
$$ 
 Z : 
 = \int_X u(x) \omega ( dx) 
 = \int_X {\bf 1}_A (x) \omega ( dx) 
 = \omega ( A ) 
$$ 
 is a Poisson random variable with intensity 
 $\lambda = \sigma (A) < \infty$, and Proposition~\ref{p111} 
 yields the following corollary. 
\begin{corollary} 
 Let $F:\Omega^X \to \real$ be a bounded random variable. 
 We have 
\begin{equation} 
\label{asac} 
 E \left[ 
 F 
 Z^n 
 \right] 
 = 
 \sum_{k=0}^n 
 S(n,k) 
 \int_{A^k} 
 E\left[ 
 \varepsilon^+_{\eufrak{s}_k } 
 F 
 \right] 
 \sigma ( ds_1 ) \cdots \sigma ( ds_k ) 
, 
 \qquad n \in \inte. 
\end{equation} 
 where $S(n,k)$ 
 denotes the Stirling number of the 
 second kind, i.e. the number of ways to partition a set of $n$ objects 
 into $k$ non-empty subsets. 
\end{corollary} 
\begin{Proof} 
 By Proposition~\ref{p111} we have 
\begin{eqnarray*} 
 E \left[ 
 F 
 Z^n 
 \right] 
 & = & 
 \sum_{ P_1,\ldots,P_k} 
 E\left[ 
 \int_{A^k} 
 \varepsilon^+_{ \eufrak{s}_k } 
 F 
 \hskip0.05cm 
 \sigma ( ds_1 ) 
 \cdots 
 \sigma ( ds_k ) 
 \right] 
\\ 
 & = & 
 \sum_{k=0}^n 
 \sum_{ 
 \substack{ 
 l_1 + \cdots + l_k = n 
 \\ 
 l_1,\ldots ,l_k \geq 1 
 } 
 } 
 {\cal N}_{\eufrak{L}_k} 
 \int_{A^k} 
 E\left[ 
 \varepsilon^+_{\eufrak{s}_k } 
 F 
 \right] 
 \sigma ( ds_1 ) \cdots \sigma ( ds_k ) 
, 
\end{eqnarray*} 
 and it remains to note that 
$$ 
 S(n,k) = 
 \sum_{ 
 \substack{ 
 l_1 + \cdots + l_k = n 
 \\ 
 l_1,\ldots ,l_k \geq 1 
 } 
 } 
 {\cal N}_{l_1,\ldots ,l_k}, 
 \qquad 
 0 \leq k \leq n. 
$$ 
\end{Proof} 
 As a consequence of Relation~\eqref{asac} 
 we find 
$$ 
 \Cov ( 
 F 
 , 
 Z^n 
 ) 
 = 
 \sum_{k=0}^n 
 S(n,k) 
 \int_{A^k} 
 E\left[ 
 ( 
 \varepsilon^+_{\eufrak{s}_k } 
 F 
 - F ) 
 \right] 
 \sigma ( ds_1 ) \cdots \sigma ( ds_k ) 
, 
 \qquad n \in \inte, 
$$ 
 and 
\begin{equation} 
\label{asdfg} 
 E [ F e^{tZ} ] 
 = 
 \sum_{k=0}^\infty 
 \frac{1}{k!} 
 ( e^t -1 )^k 
 \int_{A^k} 
 E\left[ 
 \varepsilon^+_{\eufrak{s}_k } 
 F 
 \right] 
 \sigma ( ds_1 ) \cdots \sigma ( ds_k ), 
\end{equation} 
 using e.g. Relation~(3) page 2 of \cite{sloane}. 
 Relation~\eqref{asdfg} also recovers the decomposition 
 of the Fourier transform (also called ${\cal U}$-transform) 
 on the Poisson space, cf. e.g. Proposition~3.2 of \cite{yito}. 
\\ 
 
 When $F$ has the form $F = f(Z)$ with $f:\inte \to \real$, 
 Relation~\eqref{asac} also yields an extended Chen-Stein 
 identity (see e.g. Lemma~3.3.3 of \cite{peccatitaqqu}), as 
$$ 
 E \left[ 
 Z^n 
 f(Z) 
 \right] 
 = 
 \sum_{k=0}^n 
 \lambda^k 
 S(n,k) 
 E [ f(Z+k) ] 
, 
$$ 
 and 
$$ 
 E [ f(Z) e^{tZ} ] 
 = 
 \sum_{k=0}^\infty 
 \frac{\lambda^k}{k!} 
 ( e^t -1 )^k 
 E\left[ 
 f(Z+k) 
 \right] 
. 
$$ 
 In particular, in case $F=1$, \eqref{asac} corresponds to 
 the classical relation 
$$ 
 E 
 [ 
 Z^n 
 ] 
 = 
 B_n (\lambda ) 
, 
 \qquad n \in \inte, 
$$ 
 between the Poisson moments 
 and the Bell polynomials 
\begin{equation} 
\label{bn} 
 B_n ( \lambda ) = \sum_{k=0}^n 
 \lambda^k 
 S(n,k), 
 \qquad n \in \inte, 
\end{equation} 
 cf. e.g. Proposition~3.3.2 of \cite{peccatitaqqu} and 
 references therein. 
 The comparison of \eqref{**} and \eqref{bn} 
 yields the relation 
$$ 
 S(n,k) = 
 n! 
 \sum_{ 
 \substack{ 
 r_1+2r_2+\cdots + n r_n =n 
 \\ 
 r_1+r_2+\cdots + r_n =k 
 \\ 
 r_1,\ldots , r_n \geq 0
 } 
} 
 \prod_{k=1}^n \frac{1}{(k!)^{r_k}r_k!} 
 = 
 n! 
 \sum_{ 
 \substack{ 
 r_1+2r_2+\cdots + (n-k+1) r_{n-k+1} =n 
 \\ 
 r_1+r_2+\cdots + r_{n-k+1} =k 
 \\ 
 r_1,\ldots , r_n \geq 0
 } 
} 
 \prod_{k=1}^n \frac{1}{(k!)^{r_k}r_k!} 
, 
$$ 
 cf. e.g. Proposition~2.3.4 of \cite{peccatitaqqu}. 
\\ 
 
 Similarly, Proposition~\ref{p1111} applied to 
 $u (x) = {\bf 1}_A (x)$ shows that 
$$ 
 E 
 [ 
 (Z-\lambda ) ^n 
 ] 
 = 
 \sum_{k=0}^n 
 \lambda^k 
 S_2 (n,k), 
 \qquad n \in \inte, 
$$ 
 which recovers the fact that 
 the centered moments of a Poisson random 
 variable can be written 
 using the number $S_2(n,k)$ of partitions 
 of a set of size $n$ into $k$ non-singleton subsets, 
 cf. \cite{centralmoments} and Proposition~3.3.6 of 
 \cite{peccatitaqqu}. 

\footnotesize 

\def\cprime{$'$} \def\polhk#1{\setbox0=\hbox{#1}{\ooalign{\hidewidth
  \lower1.5ex\hbox{`}\hidewidth\crcr\unhbox0}}}
  \def\polhk#1{\setbox0=\hbox{#1}{\ooalign{\hidewidth
  \lower1.5ex\hbox{`}\hidewidth\crcr\unhbox0}}} \def\cprime{$'$}

\end{document}